\documentclass{amsart}
\usepackage{amssymb,amsmath,amscd}

\newtheorem{theorem}{Theorem}[section]
\newtheorem{lemma}[theorem]{Lemma}

\newtheorem{corollary}[theorem]{Corollary}

\newtheorem{remark}[theorem]{Remark}

\begin{document}

\title{RC-graphs and a generalized Littlewood-Richardson rule.}

\author{Mikhail Kogan}

\address{Department of Mathematics, Northeastern University,
Boston, MA 02115}
\email{misha@research.neu.edu}

\begin {abstract} Using a generalization of the Schensted insertion
algorithm to rc-graphs, we provide a Littlewood-Richardson rule for
multiplying certain Schubert polynomials by Schur polynomials.  
\end{abstract}

\maketitle

\section{Introduction} RC-graphs were originally introduced by Fomin and
Kirillov in \cite{fk} in connection with Yang-Baxter equation and Schubert
calculus. Combinatorial properties of rc-graphs have been later studied by
Bergeron and Billey in \cite{bb}, where rc-graphs were further applied to
Schubert calculus, in particular, the Monk's rule was proved using a
generalization of the Schensted insertion algorithm to rc-graphs. In this
paper we use this generalized algorithm to provide a Littlewood-Richardson
rule for multiplying certain Schubert polynomials by Schur polynomials.

Let us introduce some notation. An rc-graph $R$ will be a collection of
tuples $\{(i,k)|i,k\leq n\}$, which satisfy some additional properties
(see Section 2 for precise definition). Graphically every rc-graph is
given by a table of intersecting and nonintersecting strands, so that it
represents a planar history of the permutation $w_R$ (this permutation
will permute all integers, which are not greater than $n$, such that there
exists $N$ with $w(i)=i$ for every $i\leq N$.)

Define $x^R$ to be the product of $x_k$'s, with one $x_k$ for each
$(i,k)\in R$. Then we can define a Schubert polynomial of $w$ to be
$$
S_w=\sum_{w_R=w} x^R.
$$
(The standard definition of Schubert polynomials uses divided differences
operators, but it was shown in \cite{fs} and \cite{bs} that the above
formula holds.)

Let $\mu=(\mu_1,...,\mu_n)$ be a partition with $\mu_1\geq ...\geq \mu_n$.
In Section 2 we will associate to each $\mu$ a permutation $w(\mu)$. The
Schubert polynomial $S_\mu=S_{w(\mu)}$ is the Schur polynomial of
$\mu$. Since Schubert polynomials form a basis for the ring of all
polynomials, we can write
$$
S_wS_\mu=\sum_u c^u_{w,\mu} S_u,
$$
where the sum is taken over all permutations $u$. The coefficients   
$c^u_{w,\mu}$ are known to be positive and are called the generalized
Littlewood-Richardson coefficients. 

Our goal is to provide a rule for computing Littlewood-Richardson
coefficients in the case when $w$ satisfies the following property:
$$
w(i)>w(i-1)\text{ if }i\leq 0.
$$ 
(Note that this property will imply that if the permutation $w_R=w$
satisfies the above property, then $k> 0$ for each $(i,k)\in R)$.

This rule will use the generalization of the Schensted insertion algorithm
to the case of rc-graphs given in \cite{bb}. We describe this algorithm in
detail in Section 3 and denote by $R\leftarrow k$ the result of the
insertion of a number $1\leq k\leq n$ into an rc-graph $R$, and by
$R\leftarrow Y$ the result of the insertion of a Young tableau $Y$ into
an rc-graph $R$.

The key fact, which makes the generalized Littlewood-Richardson rule
possible to prove is the following lemma, which generalizes the row
bumping lemma (see \cite{f}) in the case of the classical Schensted
algorithm. We will give precise definitions of the paths of insertions in
Section 3.3. Roughly speaking these paths are the parts of rc-graphs,
which are changed during the insertion algorithms. Let us emphasize the
fact that the following lemma does not hold for the general insertion
algorithm, but works only in the special case we consider.

\setcounter{section}{3}
\setcounter{theorem}{1}

\begin{lemma}
If $x\leq y$, then the path of $x$ is weakly to the left of the path of   
$y$ in $R\leftarrow xy$.

If $x>y$ then the path of $x$ is weakly to the right of the path of $y$ in
$R\leftarrow xy$.
\end{lemma}

This lemma plays a pivotal role in the proof of the following theorem
which gives the Littlewood-Richardson rule mentioned above.

\setcounter{section}{4}
\setcounter{theorem}{0}

\begin{theorem}
Let $w$ be a permutation, which satisfies $w(i)>w(i-1)$ for each $i\leq   
0$ and let $\mu$ be any partition. Choose any rc-graph $U$ and set
$w_U=u$. Then $c^u_{w,\mu}$ is equal to the number of pairs $(R,Y)$ of an
rc-graph $R$ and a Young tableau $Y$ with $w(R)=w$ and $\mu(Y)=\mu$, such
that $R\leftarrow Y= U$.
\end{theorem}

\setcounter{section}{1}
\setcounter{theorem}{0}

\begin{remark}{\rm It is not difficult to see that the insertion algorithm
of \cite{bb} does not work in a setting more general than in Theorem
\ref{lr}. In particular, the Pieri's formula (see \cite{bs1}, \cite{p},
\cite{so}) cannot be proved using this algorithm. But a modified insertion
algorithm for rc-graphs, which proves the Pieri's formula is constructed
in \cite{kk}.} 
\end{remark}

The paper is organized as follows. In Section 2 we recall basic
definitions and properties of rc-graphs and Young tableaux. Section 3
describes the insertion algorithm together with the proof of Lemma
\ref{weakly}. Section 4 outlines the proof of Theorem \ref{lr}. Finally,
Section 5 gives the technical details needed to prove Theorem \ref{lr}.

\medskip

{\bf Acknowledgments.} I would like to thank my advisor Victor Guillemin.
The results of this paper are a part of my Ph.D. thesis written under his
supervision. I also thank Sara Billey and Allen Knutson for many helpful
discussions and suggestions.

\section {RC-graphs and Young Tableaux.}

In this section we define rc-graphs, recall basic facts about Young
tableaux and explain why rc-graphs are just generalizations of Young
tableaux. At the end of the Section we talk about Schubert and Schur
polynomials.

We start with a definition of rc-graphs (our conventions will differ from
those of \cite{fk} and \cite{bb}). Let $R$ be a finite set of pairs of
integers $R=\{(i,k)|i\leq n,k\leq n\}$ (both $i,k$ can be negative). We
will think of $R$ as a table of intersecting and nonintersecting strands.
Strands intersect for each $(i,k)\in R$ and do not intersect otherwise.
The examples are provided on Figure 1, where we have three tables of
strands $R_1=\{(2,1),(1,1),(-1,2)\}$, $R_2=\{(3,1),(2,3),(1,2)\}$ and
$R_3=\{(3,2),(2,2),(1,2),(2,3)\}$.

\begin {picture}(200,124)

%horizontal lines

\put (22,100){\line(1,0){10}}
\put (22,84){\line(1,0){10}}
\put (22,68){\line(1,0){10}}
\put (22,52){\line(1,0){10}}
\put (22,36){\line(1,0){10}}
\put (38,100){\line(1,0){10}}
\put (38,84){\line(1,0){10}}
\put (38,68){\line(1,0){10}}
\put (38,52){\line(1,0){10}}
\put (54,100){\line(1,0){10}}
\put (54,84){\line(1,0){10}}
\put (54,68){\line(1,0){10}}
\put (70,100){\line(1,0){10}}
\put (70,84){\line(1,0){10}}
\put (86,100){\line(1,0){10}}

%vertical lines

\put (35,103){\line(0,1){10}}
\put (35,87){\line(0,1){10}}
\put (35,71){\line(0,1){10}}
\put (35,55){\line(0,1){10}}
\put (35,39){\line(0,1){10}}
\put (51,103){\line(0,1){10}}
\put (51,87){\line(0,1){10}}
\put (51,71){\line(0,1){10}}
\put (51,55){\line(0,1){10}}
\put (67,103){\line(0,1){10}}
\put (67,87){\line(0,1){10}}
\put (67,71){\line(0,1){10}}
\put (83,103){\line(0,1){10}}
\put (83,87){\line(0,1){10}}
\put (99,103){\line(0,1){10}}

%intersecting strands

\put (32,84){\line(1,0){6}}
\put (48,84){\line(1,0){6}}
\put (80,100){\line(1,0){6}}
\put (35,81){\line(0,1){6}} 
\put (51,81){\line(0,1){6}}
\put (83,97){\line(0,1){6}}

%numbers

\put(14,98){2}
\put(14,82){1}
\put(14,66){0}
\put(12,50){-1}
\put(12,34){-2}

\put(-16,79){$R_1=$}

\put(32,117){2} 
\put(48,117){1} 
\put(64,117){0}
\put(78,117){-1}
\put(94,117){-2}

%nonintersecting strands

\put(95,104){\oval(8,8)[br]}
\put(79,88){\oval(8,8)[br]}
\put(63,72){\oval(8,8)[br]}
\put(47,56){\oval(8,8)[br]}
\put(31,40){\oval(8,8)[br]}

\put(63,88){\oval(8,8)[br]}
\put(71,80){\oval(8,8)[tl]} 

\put(47,72){\oval(8,8)[br]}
\put(55,64){\oval(8,8)[tl]} 

\put(31,104){\oval(8,8)[br]}
\put(39,96){\oval(8,8)[tl]}

\put(31,72){\oval(8,8)[br]}
\put(39,64){\oval(8,8)[tl]}

\put(31,56){\oval(8,8)[br]}
\put(39,48){\oval(8,8)[tl]}

\put(47,104){\oval(8,8)[br]}
\put(55,96){\oval(8,8)[tl]}

\put(63,104){\oval(8,8)[br]}
\put(71,96){\oval(8,8)[tl]}

%second RC-graph

%horizontal lines

\put (132,100){\line(1,0){10}}
\put (132,84){\line(1,0){10}}
\put (132,68){\line(1,0){10}}
\put (132,52){\line(1,0){10}}
\put (132,36){\line(1,0){10}} 
\put (148,100){\line(1,0){10}}
\put (148,84){\line(1,0){10}}
\put (148,68){\line(1,0){10}}
\put (148,52){\line(1,0){10}}
\put (164,100){\line(1,0){10}}
\put (164,84){\line(1,0){10}}
\put (164,68){\line(1,0){10}}
\put (180,100){\line(1,0){10}}
\put (180,84){\line(1,0){10}}
\put (196,100){\line(1,0){10}}

%vertical lines

\put (145,103){\line(0,1){10}}
\put (145,87){\line(0,1){10}}
\put (145,71){\line(0,1){10}} 
\put (145,55){\line(0,1){10}}
\put (145,39){\line(0,1){10}}
\put (161,103){\line(0,1){10}}
\put (161,87){\line(0,1){10}}
\put (161,71){\line(0,1){10}}
\put (161,55){\line(0,1){10}}  
\put (177,103){\line(0,1){10}}
\put (177,87){\line(0,1){10}} 
\put (177,71){\line(0,1){10}} 
\put (193,103){\line(0,1){10}}
\put (193,87){\line(0,1){10}}  
\put (209,103){\line(0,1){10}}

%numbers

\put(124,98){2}
\put(124,82){1}
\put(124,66){0}
\put(122,50){-1}
\put(122,34){-2}

\put(142,117){2}
\put(158,117){1}
\put(174,117){0}
\put(188,117){-1}
\put(204,117){-2}

\put(94,79){$R_2=$}

%intersecting strands
\put (142,68){\line(1,0){6}}
\put (145,65){\line(0,1){6}}

\put (174,84){\line(1,0){6}}
\put (177,81){\line(0,1){6}}

\put (158,100){\line(1,0){6}}
\put (161,97){\line(0,1){6}}

%nonintersecting strands
\put(205,104){\oval(8,8)[br]}  
\put(189,88){\oval(8,8)[br]}   
\put(173,72){\oval(8,8)[br]}  
\put(157,56){\oval(8,8)[br]}  
\put(141,40){\oval(8,8)[br]} 

\put(189,104){\oval(8,8)[br]}
\put(197,96){\oval(8,8)[tl]}

\put(157,72){\oval(8,8)[br]}
\put(165,64){\oval(8,8)[tl]}

\put(141,56){\oval(8,8)[br]}
\put(149,48){\oval(8,8)[tl]}

\put(141,88){\oval(8,8)[br]}
\put(149,80){\oval(8,8)[tl]}

\put(141,104){\oval(8,8)[br]}
\put(149,96){\oval(8,8)[tl]} 

\put(157,88){\oval(8,8)[br]}
\put(165,80){\oval(8,8)[tl]}

\put(173,104){\oval(8,8)[br]}
\put(181,96){\oval(8,8)[tl]}

%third RC-graph

%horizontal lines

\put (242,100){\line(1,0){10}}
\put (242,84){\line(1,0){10}}
\put (242,68){\line(1,0){10}}
\put (242,52){\line(1,0){10}}
\put (242,36){\line(1,0){10}} 
\put (258,100){\line(1,0){10}}
\put (258,84){\line(1,0){10}}
\put (258,68){\line(1,0){10}}
\put (258,52){\line(1,0){10}}
\put (274,100){\line(1,0){10}}
\put (274,84){\line(1,0){10}}
\put (274,68){\line(1,0){10}}
\put (290,100){\line(1,0){10}}
\put (290,84){\line(1,0){10}}
\put (306,100){\line(1,0){10}}

%vertical lines

\put (255,103){\line(0,1){10}}
\put (255,87){\line(0,1){10}}
\put (255,71){\line(0,1){10}} 
\put (255,55){\line(0,1){10}}
\put (255,39){\line(0,1){10}}
\put (271,103){\line(0,1){10}}
\put (271,87){\line(0,1){10}}
\put (271,71){\line(0,1){10}}
\put (271,55){\line(0,1){10}}  
\put (287,103){\line(0,1){10}}
\put (287,87){\line(0,1){10}} 
\put (287,71){\line(0,1){10}} 
\put (303,103){\line(0,1){10}}
\put (303,87){\line(0,1){10}}  
\put (319,103){\line(0,1){10}}

%numbers

\put(234,98){2}
\put(234,82){1}
\put(234,66){0}
\put(232,50){-1}
\put(232,34){-2}

\put(252,117){2}
\put(268,117){1}
\put(284,117){0}
\put(298,117){-1}
\put(314,117){-2}

\put(204,79){$R_3=$}

%nonintersecting strands
\put(315,104){\oval(8,8)[br]}  
\put(299,88){\oval(8,8)[br]}   
\put(283,72){\oval(8,8)[br]}  
\put(267,56){\oval(8,8)[br]}  
\put(251,40){\oval(8,8)[br]} 

\put(299,104){\oval(8,8)[br]} 
\put(307,96){\oval(8,8)[tl]} 

\put(267,72){\oval(8,8)[br]} 
\put(275,64){\oval(8,8)[tl]} 

\put(251,56){\oval(8,8)[br]}
\put(259,48){\oval(8,8)[tl]}

\put(251,72){\oval(8,8)[br]}
\put(259,64){\oval(8,8)[tl]} 

\put(251,104){\oval(8,8)[br]} 
\put(259,96){\oval(8,8)[tl]} 

\put(283,104){\oval(8,8)[br]} 
\put(291,96){\oval(8,8)[tl]} 

%intersecting strands
\put (252,84){\line(1,0){6}}
\put (255,81){\line(0,1){6}}

\put (268,84){\line(1,0){6}}
\put (271,81){\line(0,1){6}}

\put (284,84){\line(1,0){6}}
\put (287,81){\line(0,1){6}}

\put (268,100){\line(1,0){6}}
\put (271,97){\line(0,1){6}}

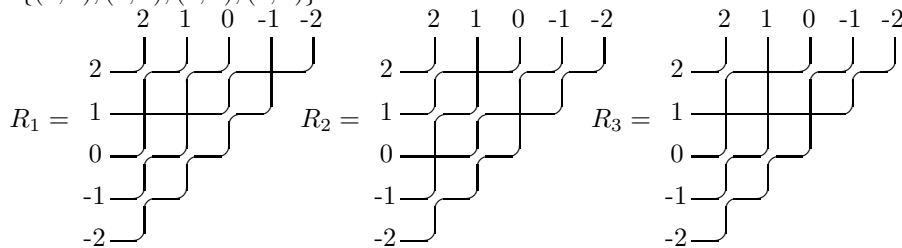
\put  (0,10){Figure 1: \small{Examples of rc-graphs.}}

\end{picture}

$R$ is called an rc-graph if no two strands intersect twice. We can think
of each rc-graph as a planar history of a permutation $w_R$, which is
defined as follows. If we label each strand by the row where it starts
from, then $w_R(i)$ is given by the column, where the $i^{\text {th}}$
strand ends. Each $w_R$ permutes all the integers, which are less than or
equal to $n$. Moreover, there always exists some negative $N$ such that
$w_R(i)=i$ for $i<N$.

For example, for the rc-graphs from Figure 1, the corresponding
permutations are given by $w_{R_1}(2,1,0,-1,-2)=(2,-2,1,0,-1)$ with
$w_{R_1}(i)=i$ for $i<-2$, $w_{R_2}(3,2,1,0,-1)=(3,1,-1,2,0)$ with
$w_{R_2}(i)=i$ for every $i<-1$ and finally
$w_{R_3}(3,2,1,0,-1)=(3,-1,1,2,0)$ with $w_{R_3}(i)=i$ for every $i<-1$.

Let us show that $R$ also provides a reduced expression for $w_R$, in
other words we can write $w_R$ as a composition of minimal number of
simple transpositions (this minimal number is called the length $\ell(w)$
of a permutation $w$). Denote by $s_i$ the simple transposition, which
permutes $i$ and $i+1$ ($i$ might be negative). Then to produce the
reduced expression for $w_R$, read each row of the rc-graph from right to
left, from the top row to the bottom one and multiply out simple
transpositions $s_{i+k-n-1}$ for each $(i,k)\in R$. It is easy to see we
get a reduced expression for $w_R$.

Let us recall some properties of rc-graphs, which were proved in 
\cite{bb}:
\begin{itemize}
\item{Let $R$ be an rc-graph, such that $(i,k),(i-1,k)$ and
$(i-1,k-\ell)$ are not in $R$ for some $i,k\leq n$ and positive $\ell$,
but all other $(i',k')$ with $k\geq k' \geq k+\ell$ and $i\geq i'\geq i+1$
are in $R$. Then we can substitute $(i,k-\ell)$ by $(i-1,k)$ without
changing $w_R$. (We also can go backwards.) These operations are called 
{\it ladder moves} of size $\ell$ at the place $(i,k)$. Examples of ladder
moves of sizes $1$ and $2$ are shown on Figure 2}
\item{For every permutation $w$ there exist a unique rc-graph $R_w$ (which
we will call a {\it top rc-graph} of $w$), such that every other
rc-graph $R$ with $w_R=w$ could be constructed from $R_w$ by a sequence
of ladder moves, which change $(i,k-\ell)$ to $(i-1,k)$ (but not the other
way).}
\item{For a top rc-graph $R_w$, if $(i,k)\in R_w$ and $i< n$, then
$(i+1,k)\in R_w$. In other words, all intersecting strands  of $R_w$ are
concentrated to the left in each row of $R_w$.}
\end{itemize}

\begin{picture}(100,100)

%ladder move of size one

\put (18,50){\line(1,0){10}}
\put (18,66){\line(1,0){10}}

\put (31,53){\line(0,1){10}}
\put (15,53){\line(0,1){10}}

\put (12,50){\line(1,0){6}}
\put (15,47){\line(0,1){6}}

\put(27,54){\oval(8,8)[br]}
\put(27,70){\oval(8,8)[br]}

\put(35,62){\oval(8,8)[tl]}
\put(19,62){\oval(8,8)[tl]}

\put (40,58){\vector(1,0){15}}
\put (55,58){\vector(-1,0){15}}

\put (68,50){\line(1,0){10}}
\put (68,66){\line(1,0){10}}

\put (81,53){\line(0,1){10}}
\put (65,53){\line(0,1){10}}

\put (78,66){\line(1,0){6}}
\put (81,63){\line(0,1){6}}

\put(77,54){\oval(8,8)[br]}
\put(61,54){\oval(8,8)[br]}

\put(69,46){\oval(8,8)[tl]}
\put(69,62){\oval(8,8)[tl]}

%ladder move of size two

\put (148,50){\line(1,0){10}}  
\put (148,66){\line(1,0){10}}   
\put (148,82){\line(1,0){10}}

\put (161,53){\line(0,1){10}}
\put (145,53){\line(0,1){10}}
\put (161,69){\line(0,1){10}}
\put (145,69){\line(0,1){10}}

\put (142,66){\line(1,0){6}}
\put (145,63){\line(0,1){6}} 
\put (158,66){\line(1,0){6}}
\put (161,63){\line(0,1){6}} 
\put (142,50){\line(1,0){6}}
\put (145,47){\line(0,1){6}} 

\put(157,54){\oval(8,8)[br]}
\put(157,86){\oval(8,8)[br]}

\put(165,78){\oval(8,8)[tl]}
\put(149,78){\oval(8,8)[tl]}

\put (170,66){\vector(1,0){15}}
\put (185,66){\vector(-1,0){15}}

\put (198,50){\line(1,0){10}}
\put (198,66){\line(1,0){10}}
\put (198,82){\line(1,0){10}}

\put (211,53){\line(0,1){10}}   
\put (195,53){\line(0,1){10}}
\put (211,69){\line(0,1){10}}
\put (195,69){\line(0,1){10}}

\put (208,82){\line(1,0){6}} 
\put (211,79){\line(0,1){6}}
\put (208,66){\line(1,0){6}}
\put (211,63){\line(0,1){6}} 
\put (192,66){\line(1,0){6}}
\put (195,63){\line(0,1){6}}

\put(207,54){\oval(8,8)[br]}
\put(191,54){\oval(8,8)[br]}

\put(199,46){\oval(8,8)[tl]}
\put(199,78){\oval(8,8)[tl]}

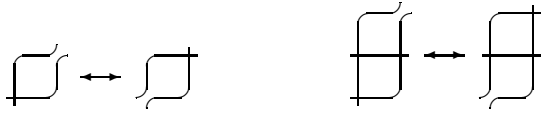
\put  (0,10){Figure 2: \small{Examples of ladder moves on rc-graphs of
sizes $1$ and $2$.}}

\end{picture}

From now on we will only work with those permutations $w$ for which
$w(i)> w(i-1)$ for each $i\leq 0$. Equivalently, every rc-graph $R$
with $w_R=w$ can be defined by the following property:
\begin{itemize}
\item{ $R$ has no two nonpositive intersecting strands.}
\end{itemize}
In particular, if $R$ satisfies the above property it lies above the
$0^{\text{th}}$ row, that is if $(i,k)\in R$ then $k\ge 1$. Let us
emphasize that starting from this point every rc-graph mentioned in this
text has to satisfy the above property. In particular, the property is
implicitly assumed in all the statements of theorems and lemmas stated
below.

Let us now define Young diagrams and tableaux. A Young diagram will be
given by a partition $\mu=(\mu_1,...,\mu_n)$, where $\mu_1\geq \mu_2 \geq
... \geq \mu_n>0$. Graphically it will be given by $\mu_i$ boxes in
$i^{\text{th}}$ row, as shown on Figure 3, where Young diagrams correspond
to partitions $(3)$, $(2,1)$ and $(3,1,1)$ respectively.

\begin{picture}(100,100)
%first diagram

\put (20,80){\line(1,0){15}}
\put (20,65){\line(1,0){15}}
\put (35,80){\line(1,0){15}}
\put (35,65){\line(1,0){15}}
\put (50,80){\line(1,0){15}}
\put (50,65){\line(1,0){15}}

\put (20,65){\line(0,1){15}}
\put (35,65){\line(0,1){15}}
\put (50,65){\line(0,1){15}}
\put (65,65){\line(0,1){15}}

% second diagram

\put (120,80){\line(1,0){15}}
\put (120,65){\line(1,0){15}}
\put (120,50){\line(1,0){15}}
\put (135,80){\line(1,0){15}}
\put (135,65){\line(1,0){15}}

\put (120,65){\line(0,1){15}}
\put (135,65){\line(0,1){15}}
\put (150,65){\line(0,1){15}}
\put (120,50){\line(0,1){15}}
\put (135,50){\line(0,1){15}}

% third diagram

\put (220,80){\line(1,0){15}} 
\put (220,65){\line(1,0){15}} 
\put (235,80){\line(1,0){15}}
\put (235,65){\line(1,0){15}}
\put (250,80){\line(1,0){15}}
\put (250,65){\line(1,0){15}} 
\put (220,50){\line(1,0){15}}
\put (220,35){\line(1,0){15}}

\put (220,65){\line(0,1){15}} 
\put (235,65){\line(0,1){15}} 
\put (250,65){\line(0,1){15}} 
\put (265,65){\line(0,1){15}}
\put (220,50){\line(0,1){15}}
\put (235,50){\line(0,1){15}}
\put (220,35){\line(0,1){15}}
\put (235,35){\line(0,1){15}}

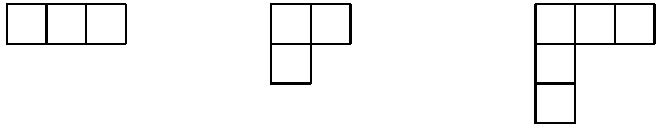
\put  (0,10){Figure 3: \small{Examples of Young diagrams.}}
\end{picture}

A Young tableaux $Y$ is a filling of a Young diagram with numbers
$1,...,n$ which satisfies the following properties. If $a$ and $b$ are two
boxes of the Young diagram, which lie in the same row, and $a$ is to the
left of $b$ then the number in $a$ is not greater than the number in $b$.
If $a$ and $b$ are in the same column and $a$ is on top of $b$, then the
number in $a$ should be strictly less than the number in $b$. We denote by
$\mu(Y)$ the partition, which corresponds to the Young tableaux $Y$.
Examples of Young tableaux are shown on Figure 4.

\begin{picture}(100,100)
%first tableaux

\put (20,80){\line(1,0){15}} 
\put (20,65){\line(1,0){15}} 
\put (35,80){\line(1,0){15}}
\put (35,65){\line(1,0){15}}
\put (50,80){\line(1,0){15}}
\put (50,65){\line(1,0){15}} 

\put (20,65){\line(0,1){15}} 
\put (35,65){\line(0,1){15}} 
\put (50,65){\line(0,1){15}} 
\put (65,65){\line(0,1){15}}

\put (25,69) {1}
\put (40,69) {1}
\put (55,69) {2}

\put (-5,69) {$Y_1=$}

% second diagram

\put (120,80){\line(1,0){15}}
\put (120,65){\line(1,0){15}}
\put (120,50){\line(1,0){15}}
\put (135,80){\line(1,0){15}}
\put (135,65){\line(1,0){15}}

\put (120,65){\line(0,1){15}}
\put (135,65){\line(0,1){15}}
\put (150,65){\line(0,1){15}}
\put (120,50){\line(0,1){15}}
\put (135,50){\line(0,1){15}}

\put (125,69) {1}
\put (140,69) {3}
\put (125,54) {3}

\put (95,69) {$Y_2=$}

%third tableaux

\put (220,80){\line(1,0){15}}
\put (220,65){\line(1,0){15}}
\put (235,80){\line(1,0){15}}
\put (235,65){\line(1,0){15}}
\put (250,80){\line(1,0){15}}
\put (250,65){\line(1,0){15}}
\put (220,50){\line(1,0){15}}
\put (220,35){\line(1,0){15}}

\put (220,65){\line(0,1){15}}
\put (235,65){\line(0,1){15}}
\put (250,65){\line(0,1){15}}
\put (265,65){\line(0,1){15}}
\put (220,50){\line(0,1){15}}
\put (235,50){\line(0,1){15}}
\put (220,35){\line(0,1){15}}
\put (235,35){\line(0,1){15}}

\put (225,69) {1}
\put (240,69) {2}
\put (255,69) {2}
\put (225,54) {2}
\put (225,39) {4}

\put (195,69) {$Y_3=$}

\put  (0,10){Figure 4: \small{Examples of Young tableux.}}
\end{picture}

Given a partition $\mu$, we construct an rc-graph $R(\mu)$ as follows. Let
$R(\mu)=\{(i,k)|1\leq k\leq n, n\geq i \geq n-\mu_k+1\}$. Then, set
$w(\mu)= w_{R(\mu)}$. (Note $R(\mu)$ is the top rc-graph of $w(\mu)$.)
Every such permutation has a unique ascent at $0$, that is $w(1)<w(0)$ but
$w(i)>w(i-1)$ if $i\neq 1$. In particular every permutation $w(\mu)$
satisfies $w(i)>w(i-1)$ if $i\leq 0$. The following lemma shows why we can
think about rc-graphs as about generalizations of Young diagrams (similar
results were obtained by Winkel in \cite{w} and pointed out in \cite{bb}).

\begin{lemma}
\label{rc=young}
RC-graphs $R$ with $w_R=w(\mu)$ are in one to one correspondence with 
Young tableaux $Y$ with $\mu(Y)=\mu$.
\end{lemma}

\proof It is easy to see that we can apply only ladder moves of size $1$
to any $R$ with $w_R=w(\mu)$. Start with the top rc-graph $R(\mu)$ and the
Young diagram, which is given by filling the $i^{\text{th}}$ row of
the Young diagram with entries $i$. Associate to each ladder move of size
$1$ an increase by $1$ of the corresponding box in the Young tableaux.
This obviously constructs a one to one correspondence between rc-graphs
with permutation $w(\mu)$ and Young tableaux with partition $\mu$.  
\endproof

Denote by $R(Y)$ the rc-graph, which is constructed out of the Young
tableau $Y$. As an illustration to the above lemma let us mention that the
first Young tableau $Y_1$ on Figure 4 correspond to the first rc-graph
$R_1$ on Figure 1. At the same time $w(\mu(Y_2))= w_{R_2}$ but $R(Y_2)
\neq R_2$.

Call any finite sequence of numbers $1,...,n$ a {\it word}. For each Young
tableau $Y$, associate a word $v(Y)$, which is given by reading the
entries of the tableau from left to right in each row, starting from the
bottom row and going to the top one. For example, the words of Young
diagrams from Figure 4 are $v(Y_1)=112$, $v(Y_2)=313$ and $v(Y_3)=42122$.

On the set of all words we define Knuth moves (originally they appeared
in \cite{k1}). These Knuth moves allow the following
changes to a word:
$$
...yxz... \Leftrightarrow  ...yzx... \ \text{ if } \ x< y \leq z
$$
and
$$
...xzy... \Leftrightarrow  ...zxy... \ \text{ if } \ x\leq y< z
$$
We say that two words $v_1$ and $v_2$ are Knuth equivalent if we can go
from one of them to another by applying a sequence of Knuth moves.

The following theorem is the key fact in the Littlewood-Richardson rule
for multiplying Schur polynomials and is very useful to us. The
proof of it can be found in \cite{f}.

\begin{theorem} 
\label{knuth words}
If $Y_1$ and $Y_2$ are two distinct Young tableaux then $v(Y_1)$ and
$v(Y_2)$ are not Knuth equivalent. Moreover, each word $v$ is Knuth
equivalent to exactly one word $v(Y)$.
\end{theorem}

Let us now recall the definitions of Schur and Schubert polynomials. Each
Young tableaux $Y$ defines a monomial $x^Y$, which is equal to the product
of $x_i$'s with one $x_i$ for each entry $i$ in the tableaux. Each
partition $\mu$ defines a Schur polynomial 
$$
S_\mu=\sum_{\mu(Y)=\mu} x^Y.
$$ 
It is well known that Schur polynomials are symmetric and that they form a
basis for the ring of symmetric polynomials in $n$ variables.

Similarly, given an rc-graph $R$ we define $x^R$ to be the product of
$x_k$'s with one $x_k$ for each $(i,k)\in R$. Then the Schubert
polynomial for the permutation $w$ is given by
$$
S_w=\sum_{w_R=w} x^R.
$$
(Let us recall again that the standard definition of rc-graphs uses
divided differences operators.) Polynomials $S_w$ form a basis for the
ring of all polynomials in $n$ variables.

Since $x^{R(Y)}=x^Y$, Proposition \ref{rc=young} implies that
$S_{w(\mu)}=S_\mu$, in other words we can think of Schubert polynomials as
generalizations of Schur polynomials.

\section{Insertion Algorithm}

The key tool in the classical Littlewood-Richardson rule for multiplying
Schur polynomials is the Schensted insertion algorithm. This algorithm was
generalized to the case of rc-graphs in \cite{bb} and used to prove Monk's
formula. We will use a special case of this generalized algorithm to
provide the Littlewood-Richardson rule for multiplying some Schubert
polynomials by Schur polynomials. This section defines the algorithm and
discusses its basic properties.

Let $R$ be an rc-graph. We would like to provide an algorithm for
inserting a number $1\leq k\leq n$ into $R$.

Let us call a pair $(i,j)$ an {\it open space}, if $(i,j)\notin R$
(two strands at position $(i,j)$ do not intersect) and the bottom strand
of the intersection is labeled by a nonpositive number, while the top 
strand is labeled by a positive number. (See Figure 5.)

\begin{picture}(100, 70)

\put (150,55){\line(1,0){10}}
\put (163,42){\line(0,1){10}}
\put (166,55){\line(1,0){10}}
\put (163,58){\line(0,1){10}}

\put(159,59){\oval(8,8)[br]} 
\put(167,51){\oval(8,8)[tl]}

\put(144,53){$a$}
\put(160,30){$b$}

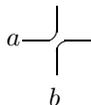
\put  (0,10){Figure 5: \small{An example of an open space: $a>0\geq b$.}}
\end{picture}

Start at the row number $k_1=k$ and find the smallest $i_1$ such that the
space $(i_1,k_1)$ is open (sometimes we will write $(i_1(k),k_1(k))$ to
indicate the dependence on $k$). Insert $(i_1,k_1)$ into $R$, in other
words make the strands intersect at $(i_1,k_1)$. If $a,b$ are the two
labels of the strands going through the place $(i_1,k_1)$ we set
$a_1(k)=a$ and $b_1(k)=b$. If we constructed an rc-graph we stop,
otherwise, it can be shown that the two stands which now intersect at the
place $(i_1,k_1)$ must also intersect at some other place $(\ell_2, k_2)$
with $k_2>k_1$. We remove $(\ell_2, k_2)$ from $R$ and find the smallest
$i_2>\ell_2$ such that $(i_2,k_2)$ is open. We insert $(i_2,k_2)$ into
$R$, set $a_2(k)$ and $b_2(k)$ to be the labels of the strands passing
through $(i_2,k_2)$ and continue the process until it stops. For
notational convenience set $k_{j+1}(k)=n+1$, if the last intersection we
inserted was $(i_j(k),k_j(k))$.

It was shown in \cite {bb} that the above algorithm stops at some point
and produces a new rc-graph, which we denote by $R\leftarrow k$. Note,
$R\leftarrow k$ and $R$ have the same number of crossings in each row,
except for the row $k$, where $R\leftarrow k$ has an additional crossing.  
Hence
$$ 
x^R x_k=x^{R\leftarrow k}. 
$$

If $v$ is a word, we denote by $R\leftarrow v$ the rc-graph we get after
inserting one by one the letters of $v$. If $Y$ is a Young tableau, we say
$R\leftarrow Y=R\leftarrow v(Y)$. Obviously we have:
$$
x^R x^Y=x^{R\leftarrow Y}.
$$

The above algorithm is a generalization of the Schensted row insertion
algorithm (see \cite{s} or \cite{f}). To prove this we just have to
translate what this algorithm means in the language of Young tableaux, in
the case when $R=R(Y)$ is constructed from some young tableau $Y$ as in
Lemma \ref{rc=young}. We omit the simple technical details of this proof,
but recall a very important fact about this algorithm (see \cite{f}). For
two Young tableaux $Y_1$ and $Y_2$ we have
\begin{equation}
\label{fact}
v(Y_1\leftarrow Y_2) \text{ is Knuth equivalent to } v(Y_1)v(Y_2).
\end{equation}
where $v(Y_1)v(Y_2)$ is just the concatenation of the two words $v(Y_1)$
and $v(Y_2)$.

Let us introduce new notations, which will be used later. During the
insertion algorithm of $k$ into $R$, each place $(i_j,k_j)$ was connected
to $(\ell_{j+1},k_{j+1})$ by two pieces of strands $s_j$ and $s^j$. (Let
us emphasize that $s_j$ and $s^j$ are just pieces of strands, which are
between rows $k_j$ and $k_{j+1}$, and their labels change during the
insertion.) We say that the left strand $s_j$ is a part of the left path
$\ell(k)$ of the insertion while the right strand $s^j$ is a part of the
right path $r(k)$ of the insertion. Both $\ell(k)$ and $r(k)$ are
collections of pieces of stands. The labeling of each piece $s_j$ in
$\ell(k)$ changes from being positive $a_j$ to nonpositive $b_j$, while
the labeling of $s^j$ in $r(k)$ change from $b_j$ to $a_j$.

Figure 6 contains an example of inserting $1$ into an rc-graph. The path
of insertion is shown on the resulting rc-graph.

\begin{picture}(100, 150)

%FIRST RC-GRAPH

\linethickness{0.4pt} 

%horizontal lines
\put (2,116){\line(1,0){10}}
\put (2,100){\line(1,0){10}}
\put (2,84){\line(1,0){10}}  
\put (2,68){\line(1,0){10}}
\put (2,52){\line(1,0){10}}
\put (2,36){\line(1,0){10}}
\put (18,116){\line(1,0){10}}
\put (18,100){\line(1,0){10}}
\put (18,84){\line(1,0){10}}  
\put (18,68){\line(1,0){10}}
\put (18,52){\line(1,0){10}}
\put (34,116){\line(1,0){10}}
\put (34,100){\line(1,0){10}}
\put (34,84){\line(1,0){10}}
\put (34,68){\line(1,0){10}}
\put (50,116){\line(1,0){10}}
\put (50,100){\line(1,0){10}}
\put (50,84){\line(1,0){10}}
\put (66,116){\line(1,0){10}}
\put (66,100){\line(1,0){10}}
\put (82,116){\line(1,0){10}}

%vertical lines

\put (15,119){\line(0,1){10}}
\put (15,103){\line(0,1){10}}
\put (15,87){\line(0,1){10}}
\put (15,71){\line(0,1){10}}
\put (15,55){\line(0,1){10}}
\put (15,39){\line(0,1){10}}
\put (31,119){\line(0,1){10}}
\put (31,103){\line(0,1){10}}
\put (31,87){\line(0,1){10}}
\put (31,71){\line(0,1){10}}
\put (31,55){\line(0,1){10}} 
\put (47,119){\line(0,1){10}}
\put (47,103){\line(0,1){10}} 
\put (47,87){\line(0,1){10}} 
\put (47,71){\line(0,1){10}} 
\put (63,119){\line(0,1){10}}
\put (63,103){\line(0,1){10}}
\put (63,87){\line(0,1){10}}
\put (79,119){\line(0,1){10}}
\put (79,103){\line(0,1){10}}
\put (95,119){\line(0,1){10}}

%numbers

\put(-6,114){4}
\put(-6,98){3}
\put(-6,82){2}
\put(-6,66){1}
\put(-6,50){0}
\put(-8,34){-1}

\put(12,135){4}
\put(28,135){3}
\put(44,135){2}
\put(60,135){1}
\put(76,135){0}
\put(90,135){-1} 

%nonintersecting strands

\put(91,120){\oval(8,8)[br]}
\put(75,104){\oval(8,8)[br]}
\put(59,88){\oval(8,8)[br]}
\put(43,72){\oval(8,8)[br]}
\put(27,56){\oval(8,8)[br]}
\put(11,40){\oval(8,8)[br]}  

\put(11,104){\oval(8,8)[br]}
\put(19,96){\oval(8,8)[tl]}

\put(11,120){\oval(8,8)[br]}
\put(19,112){\oval(8,8)[tl]}

\put(27,72){\oval(8,8)[br]}  
\put(35,64){\oval(8,8)[tl]} 

\put(11,56){\oval(8,8)[br]}  
\put(19,48){\oval(8,8)[tl]}

\put(27,104){\oval(8,8)[br]}
\put(35,96){\oval(8,8)[tl]}  

\put(27,120){\oval(8,8)[br]}
\put(35,112){\oval(8,8)[tl]} 

\put(43,104){\oval(8,8)[br]}
\put(51,96){\oval(8,8)[tl]} 

\put(43,120){\oval(8,8)[br]} 
\put(51,112){\oval(8,8)[tl]} 

\put(59,104){\oval(8,8)[br]}
\put(67,96){\oval(8,8)[tl]} 

\put(75,120){\oval(8,8)[br]}
\put(83,112){\oval(8,8)[tl]}

%intersecting strands

\put (12,68){\line(1,0){6}}
\put (15,65){\line(0,1){6}}

\put (12,84){\line(1,0){6}}
\put (15,81){\line(0,1){6}}

\put (28,84){\line(1,0){6}}
\put (31,81){\line(0,1){6}}

\put (44,84){\line(1,0){6}}
\put (47,81){\line(0,1){6}}

\put (60,116){\line(1,0){6}}
\put (63,113){\line(0,1){6}}

\put (80,68){\vector(-1,0){25}}

%SECOND RC-GRAPH

\put (122,116){\line(1,0){10}} 
\put (122,100){\line(1,0){10}} 
\put (122,84){\line(1,0){10}}  
\put (122,68){\line(1,0){10}}  
\put (122,52){\line(1,0){10}}  
\put (122,36){\line(1,0){10}}   
\put (138,116){\line(1,0){10}}
\put (138,100){\line(1,0){10}}
\put (138,84){\line(1,0){10}} 
\put (138,68){\line(1,0){10}} 
\put (138,52){\line(1,0){10}}
\put (154,116){\line(1,0){10}}
\linethickness{1pt} 
\put (154,100){\line(1,0){10}}
\linethickness{0.4pt} 
\put (154,84){\line(1,0){10}}
\linethickness{1pt} 
\put (154,68){\line(1,0){10}}
\linethickness{0.4pt} 
\linethickness{1pt} 
\put (170,116){\line(1,0){10}}
\put (170,100){\line(1,0){10}}
\linethickness{0.4pt} 
\put (170,84){\line(1,0){10}}
\put (186,116){\line(1,0){10}}
\put (186,100){\line(1,0){10}}
\put (202,116){\line(1,0){10}}

%vertical lines

\put (135,119){\line(0,1){10}}
\put (135,103){\line(0,1){10}}
\put (135,87){\line(0,1){10}}
\put (135,71){\line(0,1){10}}
\put (135,55){\line(0,1){10}}
\put (135,39){\line(0,1){10}}
\put (151,119){\line(0,1){10}}
\put (151,103){\line(0,1){10}}
\linethickness{1pt} 
\put (151,87){\line(0,1){10}}
\put (151,71){\line(0,1){10}}
\linethickness{0.4pt} 
\put (151,55){\line(0,1){10}}
\put (167,119){\line(0,1){10}}
\linethickness{1pt} 
\put (167,103){\line(0,1){10}}
\put (167,87){\line(0,1){10}}
\put (167,71){\line(0,1){10}}
\linethickness{0.4pt} 
\put (183,119){\line(0,1){10}}
\linethickness{1pt} 
\put (183,103){\line(0,1){10}}
\linethickness{0.4pt} 
\put (183,87){\line(0,1){10}}
\put (199,119){\line(0,1){10}}
\put (199,103){\line(0,1){10}}
\put (215,119){\line(0,1){10}}

%numbers

\put(114,114){4}
\put(114,98){3}
\put(114,82){2}
\put(114,66){1}
\put(114,50){0}
\put(112,34){-1}

\put(132,135){4}
\put(148,135){3}
\put(164,135){2}
\put(180,135){1}
\put(196,135){0}
\put(210,135){-1}

%nonintersecting strands

\put(211,120){\oval(8,8)[br]}
\put(195,104){\oval(8,8)[br]}
\put(179,88){\oval(8,8)[br]}
\put(163,72){\oval(8,8)[br]}
\put(147,56){\oval(8,8)[br]}
\put(131,40){\oval(8,8)[br]}

\put(131,104){\oval(8,8)[br]}
\put(139,96){\oval(8,8)[tl]}

\put(131,120){\oval(8,8)[br]}
\put(139,112){\oval(8,8)[tl]}

\put(131,56){\oval(8,8)[br]} 
\put(139,48){\oval(8,8)[tl]} 

\put(147,104){\oval(8,8)[br]}
\put(155,96){\oval(8,8)[tl]}

\put(147,120){\oval(8,8)[br]}
\put(155,112){\oval(8,8)[tl]}

\put(163,104){\oval(8,8)[br]}
\put(171,96){\oval(8,8)[tl]} 

\put(163,120){\oval(8,8)[br]}
\put(171,112){\oval(8,8)[tl]}

\put(179,104){\oval(8,8)[br]}  
\put(187,96){\oval(8,8)[tl]}   

\put(195,120){\oval(8,8)[br]}  
\put(203,112){\oval(8,8)[tl]}  

%path ovals

\put(163,72){\oval(7.3,7.3)[br]}
\put(163,72){\oval(8.7,8.7)[br]}

\put(155,96){\oval(7.3,7.3)[tl]}  
\put(163,104){\oval(7.3,7.3)[br]}
\put(171,96){\oval(7.3,7.3)[tl]}
\put(171,112){\oval(7.3,7.3)[tl]}
\put(179,104){\oval(7.3,7.3)[br]}

\put(155,96){\oval(8.7,8.7)[tl]}
\put(163,104){\oval(8.7,8.7)[br]}
\put(171,96){\oval(8.7,8.7)[tl]}
\put(171,112){\oval(8.7,8.7)[tl]}
\put(179,104){\oval(8.7,8.7)[br]}

%intersecting strands

\put (132,68){\line(1,0){6}} 
\put (135,65){\line(0,1){6}} 

\put (132,84){\line(1,0){6}} 
\put (135,81){\line(0,1){6}}

\linethickness{1pt} 
\put (148,68){\line(1,0){6}}  
\put (151,65){\line(0,1){6}}  
\linethickness{0.4pt} 

\put (148,84){\line(1,0){6}}
\linethickness{1pt} 
\put (151,81){\line(0,1){6}}
\linethickness{0.4pt} 

\put (164,84){\line(1,0){6}}
\linethickness{1pt} 
\put (167,81){\line(0,1){6}}
\linethickness{0.4pt} 

\put (180,116){\line(1,0){6}}
\put (183,113){\line(0,1){6}}

%circle

\put(151,68){\oval(9,9)}

%%THIRD RC-GRAPH

\put (242,116){\line(1,0){10}}
\put (242,100){\line(1,0){10}} 
\put (242,84){\line(1,0){10}}  
\put (242,68){\line(1,0){10}}  
\put (242,52){\line(1,0){10}}  
\put (242,36){\line(1,0){10}}  
\put (258,116){\line(1,0){10}}
\put (258,100){\line(1,0){10}}
\put (258,84){\line(1,0){10}}
\put (258,68){\line(1,0){10}}
\put (258,52){\line(1,0){10}}
\linethickness{1pt}
\put (274,116){\line(1,0){10}}
\put (274,100){\line(1,0){10}}
\linethickness{0.4pt}
\put (274,84){\line(1,0){10}}
\linethickness{1pt}
\put (274,68){\line(1,0){10}} 
\linethickness{0.4pt}
\linethickness{1pt}
\put (290,116){\line(1,0){10}}
\put (290,100){\line(1,0){10}}
\linethickness{0.4pt}
\put (290,84){\line(1,0){10}} 
\put (306,116){\line(1,0){10}}
\put (306,100){\line(1,0){10}}
\put (322,116){\line(1,0){10}}

%vertical lines

\put (255,119){\line(0,1){10}}
\put (255,103){\line(0,1){10}} 
\put (255,87){\line(0,1){10}} 
\put (255,71){\line(0,1){10}} 
\put (255,55){\line(0,1){10}}
\put (255,39){\line(0,1){10}} 
\linethickness{1pt}
\put (271,119){\line(0,1){10}}
\linethickness{0.4pt}
\put (271,103){\line(0,1){10}}
\linethickness{1pt}
\put (271,87){\line(0,1){10}}
\put (271,71){\line(0,1){10}}
\linethickness{0.4pt}
\put (271,55){\line(0,1){10}} 
\linethickness{1pt}
\put (287,119){\line(0,1){10}}
\put (287,103){\line(0,1){10}}
\put (287,87){\line(0,1){10}}
\put (287,71){\line(0,1){10}} 
\linethickness{0.4pt}
\put (303,119){\line(0,1){10}}
\linethickness{1pt}
\put (303,103){\line(0,1){10}}
\linethickness{0.4pt}
\put (303,87){\line(0,1){10}}
\put (319,119){\line(0,1){10}}
\put (319,103){\line(0,1){10}}
\put (335,119){\line(0,1){10}}

%numbers

\put(234,114){4}
\put(234,98){3}
\put(234,82){2}
\put(234,66){1}
\put(234,50){0}
\put(232,34){-1}

\put(252,135){4}
\put(268,135){3}
\put(284,135){2}
\put(300,135){1}
\put(316,135){0}
\put(330,135){-1}

%nonintersecting strands

\put(331,120){\oval(8,8)[br]}
\put(315,104){\oval(8,8)[br]}
\put(299,88){\oval(8,8)[br]}
\put(283,72){\oval(8,8)[br]}
\put(267,56){\oval(8,8)[br]}
\put(251,40){\oval(8,8)[br]}

\put(251,104){\oval(8,8)[br]}
\put(259,96){\oval(8,8)[tl]}

\put(251,120){\oval(8,8)[br]}
\put(259,112){\oval(8,8)[tl]}

\put(251,56){\oval(8,8)[br]} 
\put(259,48){\oval(8,8)[tl]} 

\put(267,104){\oval(8,8)[br]}
\put(275,96){\oval(8,8)[tl]}

\put(283,120){\oval(8,8)[br]}
\put(291,112){\oval(8,8)[tl]}

\put(283,104){\oval(8,8)[br]}
\put(291,96){\oval(8,8)[tl]} 

\put(299,104){\oval(8,8)[br]}  
\put(307,96){\oval(8,8)[tl]}   

\put(299,120){\oval(8,8)[br]}
\put(307,112){\oval(8,8)[tl]}

\put(315,120){\oval(8,8)[br]}  
\put(323,112){\oval(8,8)[tl]}  

%path ovals

\put(283,72){\oval(7.3,7.3)[br]}
\put(283,72){\oval(8.7,8.7)[br]}

\put(283,120){\oval(8.7,8.7)[br]}
\put(291,112){\oval(8.7,8.7)[tl]}
\put(283,120){\oval(7.3,7.3)[br]}
\put(291,112){\oval(7.,7.3)[tl]}

\put(275,96){\oval(7.3,7.3)[tl]}
\put(283,104){\oval(7.3,7.3)[br]}
\put(291,96){\oval(7.3,7.3)[tl]}
\put(299,104){\oval(7.3,7.3)[br]}

\put(275,96){\oval(8.7,8.7)[tl]} 
\put(283,104){\oval(8.7,8.7)[br]}
\put(291,96){\oval(8.7,8.7)[tl]} 
\put(299,104){\oval(8.7,8.7)[br]}

%intersecting strands

\put (252,68){\line(1,0){6}}   
\put (255,65){\line(0,1){6}}   

\put (252,84){\line(1,0){6}}
\put (255,81){\line(0,1){6}}

\linethickness{1pt}
\put (268,68){\line(1,0){6}} 
\put (271,65){\line(0,1){6}} 
\linethickness{0.4pt}

\put (268,84){\line(1,0){6}}
\linethickness{1pt}
\put (271,81){\line(0,1){6}}
\linethickness{0.4pt}

\put (284,84){\line(1,0){6}}
\linethickness{1pt}
\put (287,81){\line(0,1){6}}
\linethickness{0.4pt}

\linethickness{1pt}
\put (268,116){\line(1,0){6}}
\put (271,113){\line(0,1){6}}
\linethickness{0.4pt}

%circle

\put (271,116){\oval(9,9)}

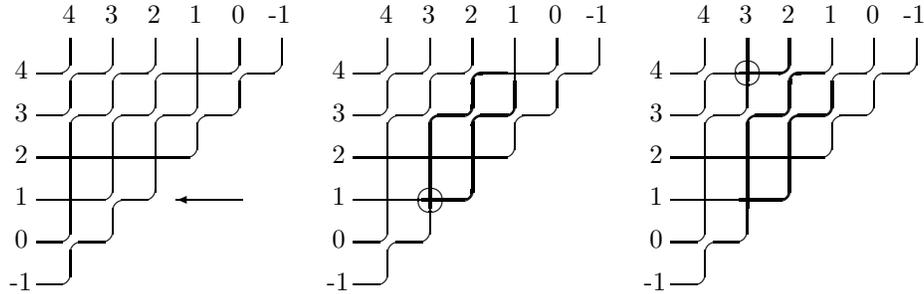
\put  (0,10){Figure 6: \small{Insertion algorithm and path of insertion.}}

\end{picture}

Recall that we are considering only those rc-graphs, for which no two
nonpositive strands intersect. The following lemma shows that the
insertion algorithm preserves this property.

\begin{lemma}
If no two nonpositive strands intersect in $R$, then no two nonpositive
strands intersect in $R\leftarrow k$.
\end{lemma}

\proof During the insertion algorithm the only possibility for introducing
new intersections of nonpositive strands is when a strand $s_j$ 
from $\ell(k)$ becomes nonpositive. 

Let us show by contradiction that $s_j$ cannot intersect any nonpositive
strand. Assume that $s_j$ (labeled by $b_j\leq 0$) in $R\leftarrow k$ is
intersected by some nonpositive strand $s$ in row $k_j< k'\leq k_{j+1}$.
Look at the whole strand $s'$ in $R\leftarrow k$, which is labeled by
$a_j$. $s'$ starts above zero, while $s$ starts below zero, at the same
time $s$ is to the left of $s'$ in the row $k'$, hence these two strands
must intersect in $R\leftarrow k$ below the row $k'$. The strand $s'$
below the row $k'$ consists of two parts: one of them is $s^j$, which was
labeled by $b_j$ in $R$, and the other one is the rest of the strand below
the row $k_j$, which is labeled by $a_j$ in both $R$ and $R\leftarrow k$.

$s$ cannot intersect $s^j$, since no two nonpositive strands intersect in
$R$. At the same time $s$ cannot intersect the rest of $s'$, since it
already intersects the strand labeled by $a_j$ in $R$ once at row $k'$,
and cannot intersect it for the second time. So we found a contradiction
and this lemma is proved.
\endproof

The lemma immediately leads to the following property
\begin{itemize}
\item{All strands passing between $r(k)$ and $\ell(k)$ are positive.}
\end{itemize}
Indeed, each strand between $r(k)$ and $\ell(k)$ has to cross at least one
strand from $r(k)$ or $\ell(k)$, but, since $r(k)$ is nonpositive in $R$
and $\ell(k)$ is nonpositive in $R\leftarrow k$, and no two nonpositive
strand can intersect, the above property holds.

Here is a very important lemma, which does not hold if we do not assume
that no two nonpositive strands intersect in $R$ (see \cite{bb} for a
counterexample).

\begin{lemma} 
\label{weakly}
If $x\leq y$, then the path of $x$ is weakly to the left of the path of
$y$ in $R\leftarrow xy$.

If $x>y$ then the path of $x$ is weakly to the right of the path of $y$ in
$R\leftarrow xy$.
\end{lemma}

\begin{remark} {\rm When we say that the path of $x$ is to the left
(right) of the path of $y$, we imply that the right path of $x$ is to the
left of the left path of $y$ (respectively, the left path of $x$ is to 
right of the right path of $y$). The word {\it weakly} stands for the
fact that $r(x)$ and $\ell(y)$ (respectively $\ell(x)$ and $r(y)$) might
have some common parts.} 
\end{remark}

\proof For the case $x\leq y$ the right path $r(x)$ of $x$ in $R\leftarrow
x$ contains strands which are all greater than zero after the insertion.
Thus when we start inserting $y$ into $R\leftarrow x$ each row $k\geq x$
should contain an open space to the right of the right path of $x$ (since
the right path of $x$ is positive). Hence the left path of $y$ is going to
stay strictly to the right of the right path of $x$, until at some point
it might happen that left path of $y$ is the same as the right path of
$x$. 

It can occur only when an open space $(i,k)=(i_j(y),k_j(y))$ in
$R\leftarrow x$ contains strands $s_j(y)$ and $s^j(y)$, such that part of
$s_j(y)$ is a part of $r(x)$. In other words, $s_j(y)$ and $s^{j'}(x)$
have a common part. If the insertion algorithm stops at this point
there is nothing to prove. Otherwise, there should be a
place $(i_{j+1}(y),k_{j+1}(y))$ where strands $s_j(y)$ and $s^j(y)$
intersect again. We would like to show 
\begin{equation}
\label{weakly1}
k_{j+1}(y) \geq k_{j'+1}(x).
\end{equation} 
This will be enough to prove the first part of the lemma. Indeed if $r(x)$
and $\ell(y)$ coincide at the row $k$, they have to separate at the row
$k_{j'+1}(x)$ by (\ref{weakly1}), so that $r(x)$ moves to the left of
$\ell(y)$. If they coincide again at some higher row, we can repeat the
argument and show that they have to separate again.

To prove (\ref{weakly1}), note $s^{j'}(x)$ was nonpositive in $R$, hence
it cannot intersect $s^j$, which was also nonpositive in $R$. Thus if
$s_j(y)$ and $s^j(y)$ intersect in $R$, it should happen above the row
where $s^{j'}(x)$ ends, in other words, above the row of intersection of
$s_{j'}(x)$ and $s^{j'}(x)$, but this row is exactly
$k_{j'+1}(x)$. Therefore (\ref{weakly1}) holds and the
first part of the lemma is proved.

In the case $x>y$, the right path of $y$ gets changed from being a set of
nonpositive strands to positive strands. The left path of $x$ in
$R\leftarrow x$ contains only nonpositive strands, so these two paths
cannot intersect (but some parts of them can coincide), since no two
nonpositive strands can intersect.

Let's now argue by contradiction that $r(y)$ is weakly to the left of
$\ell(x)$ using the fact that $r(y)$ cannot intersect $\ell(x)$. Pick the
smallest $k$, such that the right path of $y$ is to the right of the left
path of $x$. This could not happen because of an intersection of $r(y)$
and $\ell(x)$. Thus in the row $k$ the insertion of $x$ into $R$ we had to
remove some $(\ell_j(x),k_j(x))= (\ell_j(x),k)$ from $R$ and add some
$(i_j(k),k)$ to $R$, moving $\ell(x)$ to the left. But only nonpositive
strands pass in row $k$ between $\ell_j(x)$ and $i_j(k)$ (otherwise, we
would get an open space there, which is impossible), hence $r(y)$ cannot
pass between $\ell_j(x)$ and $i_j(x)$ and it must coincide with $\ell(x)$
in the row $k-1$. Moreover, the strand passing the row $k$ directly to
the left of $(\ell_j(x),k)$ is nonpositive in $R\leftarrow x$. At the same
time, the strands between right and left paths of $y$ are always positive,
so the strand passing the row $k$ directly to the left of $(\ell_j(x),k)$
is positive in $R\leftarrow x$. We found a contradiction, which means that
the second part of the lemma is proved.
\endproof

This Lemma immediately proves that if $0<x< y \leq z\leq n$ then
\begin {equation}
\label {knuth}
R\leftarrow yxz= R\leftarrow yzx
\end{equation}

Indeed in $R\leftarrow yx$ we know that $\ell(y)$ is weakly to the right 
of $r(x)$, so $r(y)$ is unchanged when we insert $x$ into $R\leftarrow
y$. At the same time when we insert $z$ in $R\leftarrow y$ the left path 
$\ell(z)$ is weakly to the right of $r(y)$. So, paths of $x$ and $z$ are
separated by the path of $y$ and, in particular, do not have any common
strands. Hence multiplication of $R\leftarrow y$ by $x$ commutes with
multiplication by $z$, which proves (\ref{knuth}).

Thus if $v_1$ and $v_2$ are two Knuth equivalent words, which can be
gotten from one another using only Knuth moves of the first type, we have
$$
R\leftarrow v_1 = R \leftarrow v_2
$$

Let us talk about how the permutation of $R$ changes after the
insertion. Notice that after each step of the algorithm the permutation
$w_R$ does not change except for the last step. At the end we make two
nonintersecting stands labeled by $c$ and $d$ intersect, which means that
$$
w_{R\leftarrow x} = w_R s_{c,d}
$$
where $s_{c,d}$ is the transposition (with $c>0\geq d$), which
interchanges the elements in position $c$ and $d$, when it acts on a
permutation from the right. Moreover,
$$
l(w_R s_{c,d})= l(w_R)+1
$$
Conversely, given an rc-graph $R'$ with $w_{R'}= w_R s_{c,d}$, such that 
$l(w_R s_{c,d})= l(w_R)+1$ and $c>0\geq d$ we can traverse the above
algorithm backwards starting by finding the unique intersection of
strands labeled by $c$ and $d$, making them nonintersecting and then
proceeding in the opposite order. For more details about the inverse of the
insertion algorithm see \cite{bb}, where Monk's formula was proved using
this inverse insertion algorithm.

\section{Littlewood-Richardson rule for multiplication Schubert
polynomials by Schur polynomials.}

Given a Schur polynomial $S_\mu$ and a Schubert polynomials $S_w$ their
product can be uniquely written as a sum of Schubert polynomials:
$$
S_w S_\mu=\sum_{u} c^u_{w,\mu} S_u,
$$
where the sum is taken over all the permutations $u$. The coefficients
$c^u_{w,\mu}$ are called Littlewood-Richardson coefficients and are known
to be positive. The following theorem provides a rule for computing these
coefficients:

\begin{theorem}
\label{lr}
Let $w$ be a permutation, which satisfies $w(i)>w(i-1)$ for each $i\leq
0$ and let $\mu$ be any partition. Choose any rc-graph $U$ and set
$w_U=u$. Then $c^u_{w,\mu}$ is equal to the number of pairs $(R,Y)$ of an
rc-graph $R$ and a Young tableau $Y$ with $w(R)=w$ and $\mu(Y)=\mu$, such
that $R\leftarrow Y= U$.
\end{theorem}

The next three lemmas will lead to the proof of the above theorem. Let us
define a Young diagram $\nu_m$ to be just one row of $m$ boxes, so
that the corresponding partition is given by one number $m$.

\begin{lemma}
\label{lemma1}
Theorem \ref{lr} holds when $\mu=\nu_m$.
\end{lemma}

\begin{remark} {\rm The above lemma is just a special case of the Pieri
formula. Since Lemma \ref{weakly} does not hold in general, Pieri formula
was conjectured but was not proved in \cite{bb}. It was later proved by
other methods in \cite{bs1}, \cite{p}, \cite{so}. In \cite{kk} the Pieri
formula is proved using a generalization of the insertion algorithm
for rc-graphs.} \end{remark}

\begin{lemma}
\label{lemma2}
The polynomials $S_{\nu_m}$ generate the ring of symmetric
polynomials in $n$ variables. So that each symmetric polynomial $S$ can
be written as
\begin{equation}
\label{sum}
S= \sum_{(m_1,...,m_k)\in M_+} S_{\nu_{m_1}}\cdots S_{\nu_{m_k}}-
\sum_{(m_1,...,m_k)\in M_-} S_{\nu_{m_1}}\cdots S_{\nu_{m_k}}
\end{equation}
where $M_+$ and $M_-$ are two sets of sequences of positive numbers.
\end{lemma}

\begin{lemma}
\label{lemma3}
Let $R$ be an rc-graph then 
$$
R\leftarrow yxz = R\leftarrow  yzx \ \text{ if }\ 0<x< y\leq z\leq n
$$
and
$$
R\leftarrow xzy =R\leftarrow  zxy\ \text{ if }\ 0<x\leq y< z\leq n.
$$
\end{lemma}

\begin{corollary} 
\label {corollary}
If $R$ is an rc-graph and $Y_1$ and $Y_2$ are two Young tableaux then
$$
R\leftarrow (Y_1\leftarrow Y_2) = (R\leftarrow Y_1)\leftarrow Y_2.
$$
\end{corollary}

We postpone the proofs of the above three lemmas until the next
Section. Let us just note that Corollary \ref{corollary} follows easily
from Lemma \ref {lemma3} and Fact (\ref{fact}).

Let us show how Theorem \ref {lr} can be proved using the above three
lemmas. We define the sets $\mathcal R_w$ and $\mathcal Y_\mu$ to be
$$
\mathcal R_w= \bigcup_{w_R=w} R \ \text{ and } \
\mathcal Y_\mu=\bigcup_{\mu(Y)=\mu} Y.
$$
If
$$
\mathcal R =\mathcal R_w \cdot \mathcal Y_\mu = \bigcup_{R\in
\mathcal R_w, Y\in \mathcal Y_\mu} R\leftarrow Y,
$$
we would like to show that
\begin {equation}
\label{equation2}
\mathcal R=\bigcup_u \mathcal R_u.
\end{equation}
This implies that each $\mathcal R_u$ is taken $c^u_{w,\mu}$ times in the
above union, since there is a unique way of writing $S_w S_\mu$ as a sum
of Schubert polynomials. Hence (\ref{equation2}) will prove the theorem.

Use Lemma \ref{lemma2} to write 
$$
S_\mu=\sum_{(m_1,...,m_k)\in M_+} S_{\nu_{m_1}}\cdots S_{\nu_{m_k}}-
\sum_{(m_1,...,m_k)\in M_-} S_{\nu_{m_1}}\cdots S_{\nu_{m_k}}
$$
this immediately implies that 
\begin{equation}
\label{equation1}
\mathcal Y_\mu =
\bigcup_{(m_1,...,m_k)\in M_+} \mathcal Y_{\nu_{m_1}}\cdots
\mathcal Y_{\nu_{m_k}} -
\bigcup_{(m_1,...,m_k)\in M_-} \mathcal Y_{\nu_{m_1}}\cdots
\mathcal Y_{\nu_{m_k}}
\end{equation}
where the minus stands for the set theoretic difference of the two sets
and where $\mathcal Y_{\mu_1}\cdot \mathcal Y_{\mu_2}
=\bigcup_{\mu(Y_1)=\mu_1, \mu(Y_2)=\mu_2} Y_1\leftarrow Y_2$. 

The reason why we can take the set theoretic difference in the above
formula is the following. By Lemma \ref{lemma1} both first and second sets
in (\ref{equation1}) could be broken up into unions of $\mathcal
Y_{\mu'}$ (since any insertion into a Young tableaux produces a Young 
tableaux). But since $S_\mu$ cannot be written as a nontrivial linear
expression of $S_{\mu'}$'s the set theoretical difference above is well
defined.

Thus we can conclude:
$$
\mathcal R =
\mathcal R_w \cdot(\bigcup_{(m_1,...,m_k)\in M_+} \mathcal
Y_{\nu_{m_1}}\cdots \mathcal Y_{\nu_{m_k}}) -
\mathcal R_w \cdot(\bigcup_{(m_1,...,m_k)\in M_-} \mathcal
Y_{\nu_{m_1}}\cdots \mathcal Y_{\nu_{m_k}})
$$
Using Corollary \ref{corollary} we can immediately see that the set
theoretic difference is well-defined in the above formula. On the other
hand, this formula and Lemma \ref{lemma1} shows that $\mathcal R$ can be
written in the form (\ref{equation2}), since by Lemma \ref{lemma1} each
$\mathcal R_w \cdot \mathcal Y_{\nu_{m_1}}\cdots \mathcal Y_{\nu_{m_k}}$
is a union $\bigcup_u \mathcal R_u$. This finishes the proof of the
Theorem \ref {lr}.

\section{Technical details in the proof of Theorem \ref{lr}.}

\proof (of Lemma \ref{lemma1}) The proof of this Lemma will just be a
combination of Monk's rule and Lemma \ref{weakly}. 

Let $Y=(1\leq a_1\leq a_2\leq ...\leq a_m \leq n)$ be a filling of the
Young diagram $\nu_m$. We can easily see from the insertion algorithm that
$$
w_{R\leftarrow Y}=w_{R\leftarrow a_1...a_m}=w_Rs_{c_1,d_1}...s_{c_m,d_m}
$$
where $c_i >0 \geq d_i$, $d_1>d_2>...>d_m$ and
$l(w_R s_{c_1,d_1}...s_{c_{m'},d_{m'}})=l(w_R)+m'$ for every $m'\leq m$.

Conversely, assume we are given an rc-graph $R'$ with
$w_{R'}=w_Rs_{c_1,d_1}...s_{c_m,d_m}$ with $c_i >0 \geq d_i$,
$d_1>d_2>...>d_m$ and $l(w_Rs_{c_1,d_1}...s_{c_{m',d_{m'}}})=l(w_R)+m'$
for every $m'\leq m$. Then we can go through the inverse insertion
algorithm and delete one by one intersections of strands $c_i$ and $d_i$.
We will get $m$ numbers $a_1,...,a_m$.

Moreover, by Lemma \ref{weakly}, $a_i\leq a_{i+1}$ (otherwise we would not
have $d_i<d_{i+1}$). Thus we have even proved a slightly better version of
the Lemma:
$$
S_w S_{\nu_m}=\sum S_{ws_{c_1,d_1}...s_{c_m,d_m}}
$$
where $c_i >0 \geq d_i$, $d_1>d_2>...>d_m$ and 
$l(w_Rs_{c_1,d_1}...s_{c_{m',d_{m'}}})=l(w_R)+m'$ for every $m'\leq m$.
\endproof

\proof (of Lemma \ref{lemma2}) This is an immediate corollary of the
Jacobi-Trudi identity (see \cite{m}). We thank Sara Billey for pointing
this out to us. 
\endproof

The rest of this Section will be concerned with the proof of Lemma
\ref{lemma3}.

Recall that the first part of Lemma \ref{lemma3} followed from Lemma
\ref{weakly}. So, we just have to prove the second part of it:
\begin{equation}
\label{need to prove}
R\leftarrow xzy= R\leftarrow zxy \ \text{ for any  } R \text{ and } 0<
x\leq y <z \leq n
\end{equation}

The path of $x$ in $R\leftarrow xz$ is weakly to the left of the path of
$z$. If it is strictly to the left of the path of $z$ (in other words the
right path $r(x)$ of $x$ has no common parts with the left path $\ell(z)$
of $z$), then clearly $R\leftarrow xz=R\leftarrow zx$ and (\ref{need to
prove}) holds. An example for this situation would be $x=1$, $y=2$, $z=3$
and $R=R_3$ (the third rc-graph from Figure 1).

Hence we just have to look at the case when right path of $x$ partially
coincides with the left path of $z$. Let's assume that the bottom row
where this happens is $k$. Then by above argument, $R\leftarrow xzy=
R\leftarrow zxy$ for all rows, which are below the row $k$.

Assume that during the insertion of $x$ into $R$ an intersection
$(i_j(x),k_j(x))$ was inserted into $R$, such that $k_j(x)<k$ but
$k_{j+1}(x)>k$. Denote by $s_1$ and $s_2$ the two pieces of strands, which
connect $(i_j(x),k_j(x))$ with $(\ell_{j+1}(x),k_{j+1}(x))$. Set
$a=a_j(x)>0\geq b=b_j(x)$. So that during the insertion of $x$ into $R$,
the labeling of $s_1$ changed from $a$ to $b$, while the labeling of $s_2$
changed from $b$ to $a$.

Assume that during the insertion of $z$ into $R\leftarrow x$ we insert an
intersection of strands at the place $(i,k)=(i_{j'}(z),k_{j'}(z))$, so
that one of the strands at $(i,k)$ is $s_2$. Denote by $s_2'$ the
piece of this strand, which connects $(i,k)$ with
$(\ell_{j'+1}(z),k_{j'+1}(z))$. ($s_2$ and $s'_2$ have a common piece
between the rows $k$ and $k_{j+1}(x)$.) Take the other strand coming out
of $(i,k)$ and denote the piece of this strand, which connects $(i,k)$
with $(i_{j'+1}(z),k_{j'+1}(z))$, by $s_3$. Clearly, $a_{j'}(z)=a$ and we
set $b_{j'}(z)=c$.

Since the path of $y$ has to sit between the right path of $x$ and the
left path of $z$ below row $k$, the strand $s_2$ has to become a part of
the left path of the insertion of $y$ into $R\leftarrow xz$. Assume it
happened at some place $(i_1,k_1)=(i_{j''}(y),k_{j''}(y))$. We claim that
$(i,k)=(\ell_{j''+1}(y),k_{j''+1}(y))$. In other words, the strands, which
pass through $(i_1,k_1)$ in $R$ have to pass through $(i,k)$ in
$R$. Indeed, if this claim does not hold, then the strand labeled by
$c\leq 0$ has to pass between the left and right paths of $y$, which is
impossible. Denote by $s_3'$ the right path, which connects $(i_1,k_1)$
with $(i,k)$, so that $s_3'$ and $s_3$ are two pieces of the same strand
in $R$.

Hence during the insertion algorithm of $y$ into $R\leftarrow xz$ we had
to remove intersection $(i,k)$ and find an open space to the left of it,
call it $(\bar i,k)= (i_{j''+1}(y),k_{j''+1}(y))$.

We have two cases: 

{\bf Case 1.} $(\bar i,k)$ is to the left of the strand $s_1$.

{\bf Case 2.} $(\bar i,k)$ is to the right of the strand $s_1$.

Before going through the proofs for both cases, let us give two
examples. Case 1 happens when we take $x=y=2$, $z=3$ and $R=R_3$ from
Figure 1. For Case 2 take $n=2$ and $R=\{(2,2)\}$ then $x=y=1$ and $z=2$
will produce Case 2.

{\it Proof of Case 1.} First of all let us note that $(\bar i,k)$ is to
the left of $s_1$ if and only if the stand $s_1$ passes exactly to the
left of strand $s_2$ in the row $k$, that is there are no other strands
between $s_1$ and $s_2$ in the row $k$. Indeed, if we had other strands
between them they had to be positive in $R$ (since they lie between right
and left paths of $x$), but then $(i-1,k)$ would be an open space, so that
$\bar i=i-1$, which contradicts the fact that $(\bar i,k)$ is to the left
of the strand $s_1$. This argument also proves that $(\bar i,k)$ is to the
right of $s_1$ if and only if $\bar i =i-1$, which will be used in the
proof of the second case.

Denote by $p_1^x$ the path of the insertion of $x$ into $R$ above the row
$k$, by $p_1^z$ the path of the insertion of $z$ into $R\leftarrow x$
above the row $k$ and by $p_1^y$ the path of the insertion of $y$ into
$R\leftarrow xz$ above the row $k$. Notice that since $(\bar i, k)$ is to
the left of $s_1$ we can conclude that $p_1^z$ is weakly to the right of
$p_1^x$, while $p_1^y$ is weakly to the left of $p_1^x$.

Let us think how $R\leftarrow zxy$ looks like in this case. When we insert
$z$ into $R$, the open space $(i_{j'}(z),k_{j'}(z))$ in the row $k$ is no
longer $(i,k)$, but it is now $(i-1,k)$. Indeed, $s_1$ is labeled by $a>
0$ while $s_2$ is labeled by $b\leq 0$ in $R$, moreover, $s_1$ passes
through the space $(i-1,k)$ and together with $s_2$ creates an open space.
So we insert $(i-1,k)$ into $R$ and denote by $p_2^z$ the path of $z$ in
$R$ above the row $k$. Notice that the paths $p_2^z$ and $p_1^x$ are
identical. When we insert $x$ into $R\leftarrow z$, at the row $k$ we have
to remove $(i-1,k)=(\ell_{j+1}(x),k_{j+1}(x))$, since $s_1$ and $s_2$
intersect at $(i-1,k)$ in $R\leftarrow z$. So, we remove $(i-1,k)$ and
insert $(\bar i,k)=(i_{j+1}(x),k_{j+1}(x))$ into $R\leftarrow z$. Denote
by $p_2^x$ the path of $x$ in $R\leftarrow z$ above the row $k$. Notice
that $p_2^x$ is identical with $p_1^y$. At the same time, the path $p_2^y$
of the insertion of $y$ into $R\leftarrow zx$ above the row $k$ will be
identical with $p_1^z$.

To summarize, we have $R\leftarrow xzy=R\leftarrow zxy$ below the row $k$.
Above the row $k$ we first insert into $R$ along the path $p_2^z=p_1^x$ in
both cases. Then we insert along $p_1^z$ and along $p_1^y$ for
$R\leftarrow xzy$ and along $p_2^x$ and along $p_2^y$ for $R\leftarrow
xzy$. But since $p_1^z,p_2^y$ are weakly to the right of $p_2^z=p_1^x$
while $p_2^x,p_1^y$ are weakly to the left of $p_2^z=p_1^x$, we can apply
Lemma \ref{weakly} to show that paths $p_1^z=p_2^y$ and $p_2^x=p_1^y$ are
separated by $p_2^z=p_1^x$ and hence it does not matter along which path
above the row $k$ we insert first. This proves $R\leftarrow
xzy=R\leftarrow zxy$ above the row $k$ and finishes the proof of Case 1.

{\it Proof of Case 2.}
Case 2 is just slightly more difficult than Case 1. 

Denote by $\tilde k$ the row where the path of the insertion of $y$ into
$R\leftarrow xz$ moves to the left of $s_1$ (this has to happen above the
row $k$, but below the row $k_{j+1}(x)$). As in Case $1$, it can be shown
that in the row $\tilde k$ the strand $s_1$ has to pass directly to the
left of $s_2$. We define paths $p_i^x$, $p_i^y$ and $p_i^z$ for $i=1,2$,
which lie above the row $\tilde k$, the same way we have done it in the
first case. We also denote by $p$ the path of the insertion of $y$ between
the rows $k$ and $\tilde k$.

Then applying the same arguments as in Case 1, we can show that
$R\leftarrow xzy=R\leftarrow zxy$ below the row $k$. We can also see that
above the row $k$, $p_2^z$ will be the path $p$ between $k$ and $\tilde k$
and $p_1^x$ above $\tilde k$. Similarly to Case 1, we apply Lemma \ref
{weakly} to show that $p_2^x$ coincides with $p_1^y$ and $p_2^y$ is the
same as $p_2^z$ above the row $\tilde k$. Hence paths of insertion in both
$R\leftarrow xzy$ and $R\leftarrow zxy$ above the row $k$ are the same.
This finishes the proof of Lemma \ref{lemma3} in the second case.


\begin{thebibliography}  {GLS}

\bibitem{bb} N. Bergeron, S. Billey.
\emph {RC-graphs and Schubert polynomials.} Experiment. Math. 2 (1993),
no. 4, 257--269.

\bibitem{bs1} N. Bergeron, F. Sotille.
\emph{A Pieri-type formula for isotropic flag manifolds.}
MSRI Preprint 1998-050, arxiv:math.CO/9810025

\bibitem{bs}  S. Billey, W. Jockusch, R. Stanley.
\emph {Some combinatorial properties of Schubert polynomials.}
J. Algebraic Combin. 2 (1993), no. 4, 345--374.

\bibitem{fk} S. Fomin, A. Kirillov.
\emph{The Yang-Baxter equation, symmetric functions, and Schubert
polynomials.} Proceedings of the 5th Conference on Formal Power Series and
Algebraic Combinatorics (Florence, 1993). Discrete Math. 153 (1996),
no. 1-3, 123--143.

\bibitem{fs} S. Fomin, R. Stanley.
\emph{Schubert polynomials and the nil-Coxeter algebra.} Adv. Math. 103
(1994), no. 2, 196--207.

\bibitem{f} W. Fulton.
\emph{Young tableaux. With applications to representation theory and
geometry.} London Mathematical Society Student
Texts, 35. Cambridge University Press, Cambridge, 1997.

\bibitem{kk} M. Kogan, A. Kumar.
\emph{A proof of Pieri's formula using generalized Schensted insertion
algorithm for rc-graphs,.} preprint

\bibitem{k1} D. Knuth
\emph{Permutations, matrices, and generalized Young tableaux.}
Pacific J. Math. 34 1970 709--727.

\bibitem{m} I. G. Macdonald
\emph{Notes on Schubert Polynomials}
Montreal, Canada : Departement de mathematiques et d'informatique,
Universite du Quebec a Montreal, 1991

\bibitem{p} A. Postnikov
\emph{On a Quantum Version of Pieri's Formula}, preprint March 23, 1997.

\bibitem{s} C. Schensted.
\emph{Longest increasing and decreasing subsequences.} Canad. J. Math. 13
1961 179--191.

\bibitem{so} F. Sottile.
\emph{Pieri's formula for flag manifolds and Schubert polynomials},
Annales de l'Institut Fourier {\bf 46} (1996), 89-110.

\bibitem{w} R. Winkel.
\emph{A combinatorial bijection between Standard Young Tableaux and
reduced words of Grassmannian permutations.} 
Sem. Loth. Comb. B36h (1996).


\end{thebibliography}
\end{document}